\documentclass[11pt%,referee
]{article}

\usepackage{amsmath}

\usepackage{graphicx,color}
\usepackage{times}

\usepackage[colorlinks=true, breaklinks=true, pdfstartview=FitV, linkcolor=red, citecolor=blue, urlcolor=black]
{hyperref}

\usepackage[T1]{fontenc}

\allowdisplaybreaks

\begin{document}

\title{ A torsion-based solution to the hyperbolic regime of the $J_2$-problem
\thanks{Preprint submitted to \emph{Nonlinear Dynamics} July 23, 2022 }
}

\author{Martin Lara%
\thanks{GRUCACI, University of La Rioja, and Space Dynamics Group -- UPM}
\thanks{\tt{mlara0@gmail.com} }
\and
Alessandro Masat%
\thanks{Department of Aerospace Science and Technology, Politecnico di Milano, Milano, Italy}
\thanks{PhD Candidate, \tt{alessandro.masat@polimi.it} }
\and
Camilla Colombo%
\thanks{Department of Aerospace Science and Technology, Politecnico di Milano, Milano, Italy}
\thanks{Associate Professor, \tt{camilla.colombo@polimi.it} }
}

\date{Manuscript of \today} %, \color{red}in preparation }

\maketitle

\begin{abstract}
A popular intermediary in the theory of artificial satellites is obtained after the elimination of parallactic terms from the $J_2$-problem Hamiltonian. The resulting quasi-Keplerian system is in turn converted into the Kepler problem by a torsion. When this reduction process is applied to unbounded orbits the solution is made of Keplerian hyperbolae. For this last case, we show that the torsion-based solution provides an effective alternative to the Keplerian approximation customarily used in flyby computations. Also, we check that the extension of the torsion-based solution to higher orders of the oblateness coefficient yields the expected convergence of asymptotic solutions to the true orbit.

%\footnote{\color{red}Reviewers: Martinusi, Giorgilli, Gurfil, Healy, Alfriend, Ferrer, Maciejewski, Floria, Hao Zhang, Marco Nunes}
\end{abstract}


%\linenumbers

\section{Introduction}

The hyperbolic regime of an artificial satellite of an oblate body was solved in closed form by Hori \cite{Hori1961}, who showed in this way that the perturbation approach may be useful without restriction to the typical case of quasi-periodic motion.\footnote{As in the case of bounded motion \cite{Jezewski1983}, the invariant manifold of equatorial orbits of the main problem is integrable for unbounded motion too \cite{MartinusiGurfil2013}; see \cite{DangLuoShiZhang2019} for a general discussion.} However, at difference from the classical elliptic regime, perturbed hyperbolic motion is subjected to strict boundary conditions that must be imposed to the analytical approach. More precisely, the effects of the oblateness perturbation vanish when the artificial satellite approaches infinity in its hyperbolic-type motion. In consequence, the transformation between mean and osculating variables that encapsulates the perturbation solution must become the identity when evaluated at infinity.
\par

In Hori's first-order approach the $J_2$-problem Hamiltonian is directly reduced to the Kepler Hamiltonian in (mean) hyperbolic Delaunay variables. With this choice, the solution of the mean-to-osculating transformation that yields the complete Hamiltonian reduction depends unavoidably on non-periodic terms of the true anomaly. In spite of that, the explicit appearance of the true anomaly in Hori's solution is not of concern since its value is bounded by the asymptotes of a Keplerian hyperbola. Therefore, these non-periodic terms do not play any secular role. After application of the boundary conditions to Hori's transformation these kinds of terms appear always like coefficients involving the difference between the instantaneous true anomaly and the constant true anomaly of the asymptote of the Keplerian hyperbola. While non-periodic, this type of coefficient resembles the equation of the center, given by the difference between the true and mean anomalies, which comprises essential short-period terms of the elliptic motion. But the resemblance goes beyond a formal likeness, and we find that these terms can be handled advantageously using known techniques for efficiently dealing with the equation of the center.
\par

In particular, it is known that the equation of the center does not play any role in the elimination of the parallax simplification, which, on the contrary, depends only on trigonometric functions of the true anomaly, cf.~\S 6.2.1 of \cite{Lara2021}. That this last feature is not exclusive of the case of bounded perturbed Keplerian motion and also applies to unbounded motion becomes evident from Deprit's original formulation of the elimination of the parallax in polar variables \cite{Deprit1981}. Furthermore, Deprit's radial intermediary, which is obtained after the elimination of the parallax up to first order effects of the oblateness coefficient, provides an integrable solution of the $J_2$-problem, also called the \emph{main problem} in the theory of artificial satellites \cite{Brouwer1959}, that is valid for both bounded and unbounded orbits. Nevertheless, the nature of the solution is different  in the two distinct regimes, and Deprit's reduction of the main problem to a quasi-Keplerian system must be amended in the second case in order to take the boundary conditions that characterize hyperbolic-type motion into account.
\par

We compute the additional terms that are needed in the elimination of the parallax transformation for converting the hyperbolic regime of the main problem into a quasi-Keplerian system, which, in turn, is converted into a pure Keplerian system by the standard torsion transformation \cite{Deprit1981}. When the new transformation is applied to convert the mean hyperbolic Delaunay variables into osculating ones, the transformation equations are simpler than those provided by Hori. Still, the new transformation shares analogous deficiencies to those of Hori's solution. Namely, it fails in the limit case of parabolic orbits, and the accuracy of the perturbation solution deteriorates for values of the eccentricity close to 1. In this regards, it is worth mentioning that existing perturbation solutions of bounded motion suffer analogous limitations too.
\par

Moreover, we show that the reduction to a quasi-Keplerian hyberbola is easily extended to higher orders, yet the appearance of non-periodic functions of the true anomaly in the transformation equations is now unavoidable. 
While these higher order refinements of the perturbation solution may not have relevant effects in practice ---for flyby design or gravity assist maneuvers--- they radically reduce the errors of the perturbation solution close to, and after perigee passage, and serve us to illustrate how the perturbation approach follows the expected convergence beyond its customary application to the case of quasi-periodic motion. Needless to say that analogous improvements are expected if Hori's perturbation solution were eventually extended beyond the first order of the oblateness coefficient explicitly computed by himself. 
\par

Due to the hyperbolic character of unbounded perturbed Keplerian motion, one might expect high sensitivity of the main problem dynamics with respect to initial conditions. However, for common values of the $J_2$ coefficient of solar system bodies the perturbation solution succeeds in capturing these dynamics, and can even predict the qualitative features of branches of the classical Newtonian rosetta when the energy becomes positive.\footnote{In the case of zero energy, these branches are sometimes called ``fish'' orbits \cite{MartinusiGurfil2013}} In this extreme case, high orders of the perturbation solution are indeed required.
\par

\section{Deprit's radial intermediary}

The main problem of the artificial satellite admits Hamiltonian formulation. When using polar canonical variables, it is written in the form
\begin{equation} \label{Hamain}
\mathcal{M}=\frac{1}{2}\left(R^2+\frac{\Theta^2}{r^2}\right)-\frac{\mu}{r}-\frac{1}{4}J_2\frac{\mu}{r}\frac{\alpha^2}{r^2}\left(2-3\sin^2I+3\sin^2I\cos2\theta\right),
\end{equation}
where the physical parameters $\mu$, $\alpha$, and $J_2$ characterize the gravitational field of the attracting body and denote the gravitational parameter, equatorial radius, and oblateness coefficient, respectively. The polar canonical coordinates $(r,\theta)$ stand for the radius and the polar angle, respectively, and their conjugate momenta $(R,\Theta)$ for the radial velocity and the specific angular momentum, respectively; $I=\arccos(N/\Theta)$ is the inclination of the instantaneous orbital plane, and the canonical pair $(\nu,N)$ represents the right ascension of the ascending node and the third component of the angular momentum vector per mass unit, respectively. Because $\nu$ is an ignorable coordinate of Hamiltonian (\ref{Hamain}), we immediately infer that $N$ is an integral of the main problem dynamics.
\par

Up to first order effects of $J_2$, Hamiltonian (\ref{Hamain}) can be reduced to the integrable quasi-Keplerian Hamiltonian in new, prime variables
\begin{equation} \label{Hamainp}
\mathcal{M}=\frac{1}{2}\left(R'^2+\frac{\Theta'^2}{r'^2}\right)-\frac{\mu}{r'}
-\frac{1}{2}\frac{\Theta'^2}{r'^2}\frac{1}{2}J_2\frac{\alpha^2}{(\Theta'^2/\mu)^2}\left[2-3\sin^2I(\Theta',N)\right],
\end{equation}
which is known as Deprit's radial intermediary \cite{Deprit1981}. In the case of bounded motion this intermediary has been proposed for different applications \cite{CoffeyAlfriend1984,GurfilLara2014CeMDA,HautesserresLaraCeMDA2017}. The reduction of the main problem to Deprit's radial intermediary is achieved by an infinitesimal contact transformation
$(r,\theta,\nu,R,\Theta,N,J_2)\mapsto(r',\theta',\nu',R',\Theta',N)$ of generating function
\begin{equation} \label{Gemain1}
\mathcal{U}_1=-\Theta\frac{1}{8}\frac{\alpha^2}{p^2}\left[\left(4\kappa+3\right)s^2\sin2\theta
+(4-6s^2-2s^2\cos2\theta)\sigma\right]+\mathcal{C}_0,
\end{equation}
where $s$ abbreviates the sine of the inclination, $p=\Theta^2/\mu$ is the parameter of the conic, and the non-dimensional functions
\begin{equation}
\sigma=\frac{pR}{\Theta}, \qquad \kappa=\frac{p}{r}-1,
\end{equation}
represent the projections of the eccentricity vector in the orbital frame. Namely, $\sigma=e\sin{f}$ and $\kappa=e\cos{f}$, with $e$ denoting orbital eccentricity and $f$ is the true anomaly, from which we readily obtain $\lim_{r\to\infty}\cos{f}=-1/e$.
\par

Differences between Eq.~(\ref{Gemain1}) and the equivalent term in \cite{Deprit1981} are a consequence of a distinct arrangement of the involved functions, on the one hand, and the addition of the ``integration constant'' $\mathcal{C}_0$, on the other. The latter must fulfill the condition of being free from the explicit appearance of $\theta$ when it is formulated in the specific algebra of the functions
\begin{equation} \label{statef}
C=\frac{\Theta}{p}\left(\kappa\cos\theta+\sigma\sin\theta\right),\qquad
S=\frac{\Theta}{p}\left(\kappa\sin\theta-\sigma\cos\theta\right),
\end{equation}
used in \cite{Deprit1981}. However, their use may obscure the real nature of this integration constant, which is no other than being free from short-period effects \cite{LaraSanJuanLopezOchoa2013b}. For this reason, but also for simplicity, carrying out the elimination of the parallax in Delaunay variables, rather than in the set of polar variables, is nowadays encouraged \cite{LaraSanJuanLopezOchoa2013c}.
\par

Delaunay variables are traditionally denoted $(\ell,g,h,L,G,H)$ and represent the mean anomaly, the argument of the periapsis, the longitude of the ascending node, the so-called Delaunay action, the specific angular momentum, and the third component of the angular momentum vector per mass unit, respectively. When the generating function $\mathcal{U}_1$ in Eq.~(\ref{Gemain1}) is reformulated in Delaunay variables, one must distinguish between the cases of bounded and unbounded motion. For the former the classical set of elliptic Delaunay variables is used \cite{Delaunay1860,DepritRom1970,FerrerLara2010b}, whereas hyperbolic Delaunay variables are used for the latter instead \cite{Hori1961,Floria1995}. While in both cases Eq.~(\ref{Gemain1}) takes the form
\begin{align} \nonumber
\mathcal{U}_1=& 
-G\frac{1}{8}\frac{\alpha^2}{p^2}\left\{s^2[3e\sin(f+2g)+3\sin(2f+2g)+e\sin(3f+2g)] \right. \\ \label{GemainD}
& \left. -(6s^2-4)e\sin{f}\right\}+\mathcal{C}_0(g,L,G,H),
\end{align}
the meaning of the symbols used is slightly different. Thus, while $p=G^2/\mu$ and $I=\arccos{H}/{G}$ have the same meaning in both regimes of motion, the eccentricity is $e=(1-\eta^2)^{1/2}<1$, with $\eta=G/L$, for bounded motion, whereas $e=(1+\eta^2)^{1/2}>1$ and $\eta=-G/L$ in the hyperbolic case. The true anomaly $f$ is an implicit function of the mean anomaly and the eccentricity, whose computation requires the preliminary solution of Kepler's equation. For unbounded motion, Kepler equation takes the form
\begin{equation} %\label{K:Keplerh}
\ell=e\sinh{u}-u,
\end{equation}
where the hyperbolic anomaly $u$ is related to the true anomaly by means of
\begin{equation}
\sqrt{e+1}\tanh\left(\mbox{$\frac{1}{2}$}u\right)=\sqrt{e-1}\tan\left(\mbox{$\frac{1}{2}$}f\right).
\end{equation}
\par
%{\color{magenta} no se si es un pegote y llegar\'ia con poner el Kepleriano \dots
%Then, the perturbation problem $\mathcal{M}=\mathcal{M}_{0,0}+J_2\mathcal{M}_{1,0}$ is written in hyperbolic Delaunay variables as
%\begin{align} \label{Keplerian}
%\mathcal{M}_{0,0}= & +\frac{\mu}{2a} \\
%%\mathcal{M}_{1,0}= & -\frac{1}{4}\frac{\mu}{r}\frac{\alpha^2}{r^2}\left[2-3s^2+3s^2\cos(2f+g)\right] \\
%\mathcal{M}_{1,0}= & -\frac{1}{4}\frac{\mu}{p}\frac{\alpha^2}{r^2}(1+e\cos{f})\left[2-3s^2+3s^2\cos(2f+g)\right]
%\end{align}
%where $a=L^2/\mu$ is the semi-transverse axis of the osculating hyperbola, and $r=p/(1+e\cos{f})$.
%}
%\par

\section{Polar variables mapping to the quasi-Keplerian system}

In his perturbation approach Hori chose that the new Hamiltonian be just the Keplerian part of the perturbation problem. That is, the transformed Hamiltonian in the prime hyperbolic Delaunay variables reads $\mathcal{M}'=+\mu/(2a)$, where $a={L'}^{2}/\mu$ is the semi-transverse axis of the osculating hyperbola. Because of this choice, non-trigonometric terms of the true anomaly show up in the generating function, which Hori solved up to an integration constant, as we did in Eq.~(\ref{GemainD}). At difference from the elliptic case the latter is in no way arbitrary, and Hori determined it by imposing the required boundary conditions at infinity to his transformation of the hyperbolic Delaunay variables.% After that, he computed the first order transformation equations.
\footnote{Note a typo in Hori's final transformation equation for the longitude of the node, where the exponent of $G'$ in the first coefficient must be 4 instead of 3, as readily shown by checking dimensions.}
\par

Hori had not yet invented his perturbation method based on Lie series \cite{Hori1966} and relied on von Zeipel's algorithm for computing his first order solution \cite{vonZeipel1916}. Rather, we resort to the method of Lie transforms \cite{Deprit1969} in following derivations. Moreover, in contrast with Hori's choice of the Keplerian, we are satisfied with finding the transformation that yields Deprit's radial intermediary, in this way simplifying the computation of the mapping.
\par

At the end, the quasi-Keplerian system (\ref{Hamainp}) is converted into a pure Keplerian system by means of a torsion transformation \cite{Deprit1981}. Since the torsion is carried out in polar variables, we find advantageous to derive the transformation equations of the polar variables, rather than the hyperbolic Delaunay ones, directly from Eq.~(\ref{GemainD}). That is, to the first order
\begin{equation} \label{tx1}
\xi=\xi'+J_2\xi_{0,1}|_{\xi=\xi'}
\end{equation}
where $\xi\equiv(r,\theta,\nu,R,\Theta,N)^\tau$, $\tau$ means transposition, $\xi_{0,1}\equiv\{\xi,\mathcal{U}_1\}$, and the curly brackets denote the Poisson bracket operator. Still, since the boundary conditions are naturally formulated in terms of the true anomaly, we find convenience in writing the transformation in terms of the hyperbolic Delaunay variables. To this end, we need to avail ourselves with the partial derivatives of the polar variables with respect to the hyperbolic Delaunay ones. They are listed in Table \ref{t:partials}. Note that, while mostly analogous to the case of elliptic Delaunay variables, sign changes appear in the partial derivatives with respect to the hyperbolic Delaunay action, cf.~Table 5.1 of \cite{Lara2021}.
\par

\begin{table}[htb]
\centering \small
\begin{tabular}{@{}lccc@{}}
 & $\partial/\partial\ell$ & $\partial/\partial{L}$ & $\partial/\partial{G}$ \\
\hline\noalign{\smallskip}
$\theta$ & $\displaystyle\frac{p^2}{r^2}\frac{1}{\eta^3}$
 & $\displaystyle\eta\frac{\partial\theta}{\partial{G}}$
 & $\displaystyle-\frac{1}{e^2}\frac{p R}{\Theta^2}\left(\frac{p}{r}+1\right)$ \\[2ex]
$r$ & $\displaystyle\frac{p^2 R}{\eta ^3 \Theta }$
 & $\displaystyle\eta\left(\frac{\partial{r}}{\partial{G}}-2\frac{r}{\Theta}\right)$
 & $\displaystyle\frac{1}{e^2}\frac{p}{\Theta}\left(\frac{p}{r}-1\right)$ \\[2ex]
$R$ & $\displaystyle\frac{e^2}{\eta^3}\frac{\Theta^2}{r^2}\frac{\partial{r}}{\partial{G}}$
 & $\displaystyle\eta\left(\frac{\partial{R}}{\partial{G}}+\frac{R}{\Theta}\right)$
 & $\displaystyle-\frac{1}{e^2}\frac{p^2}{r^2}\frac{R}{\Theta}$ \\
\noalign{\smallskip}\hline
\end{tabular}
\caption{Non-vanishing partial derivatives of the polar variables with respect to the hyperbolic Delaunay variables ($\partial\Theta/\partial{G}=1$). }
\label{t:partials}
\end{table}

We readily find that the correction to the radial velocity automatically fulfills the boundary condition $\{R,\mathcal{U}_1\}|_{\infty}\equiv0$. Next, we compute
\[
\{\Theta,\mathcal{U}_1\}|_{\infty}= G\frac{1}{4}\frac{\alpha^2}{p^2}\frac{s^2}{e^2}
\left[2\eta^3\sin2g-(3\eta^2+1)\cos2g\right]-\frac{\partial{\mathcal{C}_0}}{\partial{g}},
\]
from which $\mathcal{C}_0$ is solved to make the right hand of the equation to vanish. We obtain
\[
\mathcal{C}_0=-G\frac{1}{8}\frac{\alpha^2}{p^2}\frac{s^2}{e^2}\left[2\eta^3\cos2g+(3e^2-2)\sin2g\right]
+\mathcal{C}_1(L,G,H),
\]
which is the same as the function $Z$ computed by Hori except for a sign due to the different convention in the Hamiltonian formulation used by Hori. It follows the computation of
\[
\{\nu,\mathcal{U}_1\}|_{\infty}= \frac{3}{2}c\frac{\alpha^2}{p^2}\eta+\frac{\partial{\mathcal{C}_1}}{\partial{H}},
\]
where $c=H/G$. The constant $\mathcal{C}_1$ is then solved after imposing the vanishing of the right side. We obtain
\[
\mathcal{C}_1=-G\frac{3}{4}\frac{\alpha^2}{p^2}\eta(1-s^2)+\mathcal{C}_2(L,G).
\]
Analogously,
\[
\{\theta,\mathcal{U}_1\}|_{\infty}=\frac{1}{2}\eta\frac{\alpha^2}{p^2}+\frac{\partial{\mathcal{C}_2}}{\partial{G}},
\]
from which, recalling that $\eta=-G/L$ and $p=G^2/\mu$, we obtain
\[
\mathcal{C}_2=G\frac{1}{4}\frac{\alpha^2}{p^2}\eta+\mathcal{C}_3(L).
\]
Finally, we check that the choice $\mathcal{C}_3=0$ fulfills the last condition $\{r,\mathcal{U}_1\}|_{\infty}=0$. In summary,
\begin{equation} \label{c0}
\mathcal{C}_0=G\frac{1}{4}\frac{\alpha^2}{p^2}\left\{(3s^2-2)\eta-\frac{s^2}{e^2}\left[\eta^3\cos2g+\frac{1}{2}(3e^2-2)\sin2g\right]\right\},
\end{equation}
%in hyperbolic Delaunay variables, or: 
%CONFUSO. PARA QUE SE VEA CONSTANTE DEBES USAR EL ALGEBRA DE DEPRIT ....
%\begin{align} \nonumber
%\mathcal{C}_0= & G\frac{1}{4}\frac{\alpha^2}{p^2}\Big\{(3s^2-2)\eta+\frac{s^2}{e^4}\Big[
%\eta^3\left[\left(\sigma^2-\kappa^2\right)\cos2\theta-2\kappa\sigma\sin2\theta\right]
% \\ \nonumber
%& +\frac{1}{2}(3e^2-2)\left[\left(\sigma ^2-\kappa^2\right)\sin2\theta+2\kappa\sigma\cos2\theta\right)
%\Big]\Big\},
%\end{align}
thus making Eq.~(\ref{GemainD}) fully determined.
\par

It is worth noting that Eq.~(\ref{c0}) can be rewritten in polar variables using the alternative form
\[
\mathcal{C}_0=\Theta\frac{1}{4}\frac{\alpha^2}{p^2}\left\{\left(3s^2-2\right)\eta-
\frac{p^2s^2}{e^4\Theta^2}\left[\eta^3\left(C^2-S^2\right)+\left(3e^2-2\right)CS\right]
\right\},
\]
with $e=(p/\Theta)\sqrt{C^2+S^2}$, whereas $C\equiv{C}(r,\theta,R,\Theta)$ and $S\equiv{S}(r,\theta,R,\Theta)$ are given in Eq.~(\ref{statef}). In this way it is shown that $\mathcal{C}_0$ is free from the \emph{explicit} appearance of $\theta$, like it must be when using Deprit's algebra based on these functions \cite{Deprit1981}.
\par

Once the generating function has been fully determined, we compute the corrections
\begin{align} 
    \nonumber
r_{0,1} = 
& \; p\frac{1}{4}\frac{\alpha^2}{p^2}\Big\{
(3s^2-2)\left(1+\frac{e}{\eta}\sin{f}\right)
 +\frac{s^2}{2e^3}\left[ (e^2-4)\eta\sin(f-2g) \right.  \\ \nonumber
& -3e^2\eta\sin(f+2g)+(3e^2-4)\cos(f-2g)+3e^2\cos(f+2g) \\ \label{r01}
& \left. +2e^3\cos(2f+2g) \right] \Big\}
 \\ \nonumber
\theta_{0,1} = 
& \; \frac{1}{16}\frac{\alpha^2}{p^2}\Big(
\frac{1}{\eta}\Big\{6\left[2(5s^2-4)-(7s^2-6)e^2\right] +8e(3s^2-2)\cos{f} 
+2e^2 \\ \nonumber
& \times(3s^2-2)\cos2f \Big\}
+\frac{\eta}{e^3}\Big\{ (e^2-4)es^2\cos(2f-2g) +4(e^2-4)s^2 \\ \nonumber
& \times\cos(f-2g) +2 e \left[e^2(7s^2-4)-4(4 s^2-1)\right]\cos2g
-12e^2s^2 \\ \nonumber
& \times\cos(f+2g) -3e^3s^2\cos(2f+2g) \Big\}
+\frac{1}{e^3}\Big\{ (4-3e^2)es^2\sin(2f-2g) \\ \nonumber
& -4(3e^2-4)s^2\sin(f-2g) +2e\left[3e^2(5s^2-2)-4(4s^2-1)\right]\sin2g \\ \nonumber
& -8e^4(6s^2-5)\sin{f} +4e^2\left[e^2(5s^2-3)-3s^2\right]\sin(f+2g) \\ \label{z01}
& +e^3(11s^2-12)\sin(2f+2g)  +4 e^4(s^2-1)\sin(3f+2g) \Big\}
\Big)
    \\ \nonumber
\nu_{0,1} = & \; c\frac{1}{4}\frac{\alpha^2}{p^2}\Big\{
\frac{1}{e^2}\left[(3e^2-2)\sin2g+2\eta^3\cos2g\right]
-6 (\eta+e\sin{f}) \\ \label{h01}
& +3e\sin(f+2g)+3\sin(2f+2g)+e\sin(3f+2g) \Big\}
    \\ \nonumber
R_{0,1} = & \; \frac{G}{p}\frac{1}{32}\frac{\alpha^2}{p^2}\Big\{
\frac{e}{\eta}(3s^2-2)\left[2e^2\cos3f+8e\cos2f+(6e^2+8)\cos{f}+8e\right] \\ \nonumber
& +\eta\frac{s^2}{e^3}\left[
(e^2-4)e^2\cos(3f-2g) +4(e^2-4)e\cos(2f-2g)
-(e^4+4e^2 \right. \\ \nonumber
& +16)\cos(f-2g) -8(e^2+2)e\cos2g -(5e^2+16)e^2\cos(f+2g)
-12 \\ \nonumber
& \left. \times e^3\cos(2f+2g) -3e^4\cos(3f+2g) \right]
 -\frac{s^2}{e^3}\left[ (3e^2-4)e^2\sin(3f-2g) \right.  \\ \nonumber
& +4(3e^2-4)e\sin(2f-2g) +(3e^4+4e^2-16)\sin(f-2g)
+4(e^4+4) \\ \nonumber
& \times e\sin2g
 +(19e^2+16) e^2\sin(f+2g) +4(2e^2+7)e^3\sin(2f+2g) \\ \label{v01}
&  \left. +19e^4\sin(3f+2g) +4e^5\sin(4f+2g) \right] \Big\} 
\\  \nonumber
\Theta_{0,1} = & \; \Theta\frac{1}{4}\frac{\alpha^2}{p^2}s^2\Big\{
\frac{1}{e^2}\left[(3e^2-2)\cos2g-2\eta^3\sin2g\right] \\ \label{G01}
& +3e\cos(f+2g) +3\cos(2f+2g) +e\cos(3f+2g)  \Big\}.
\end{align}
The appearance of $\eta$ in denominators of $r_{0,1}$, $\theta_{0,1}$ and $R_{0,1}$ will deteriorate the efficiency of these corrections for orbits close to parabolic.
% This negative effect is more severe when using Hori's corrections to the hyperbolic Delaunay variables, in which case the correction of the Delaunay action has the factor $(a/p)$
Note that the dependence of Eqs.~(\ref{r01})--(\ref{G01}) on the true anomaly is only through trigonometric functions, in clear contrast with Hori's corrections \cite{Hori1961}. On the other hand, the explicit appearance of $f$ in Hori's solution is not criticizable as far as the possible values of the true anomaly are constrained by the values taken by the asymptotes of the Keplerian hyperbola. 
\par

%\section{Torsion}
The analytical solution is completed by applying a torsion to the quasi-Keplerian Hamiltonian (\ref{Hamainp}). Disregarding the prime notation without risk of confusion, Eq.~(\ref{Hamainp}) is rearranged like
\begin{equation} \label{Hamainq}
\mathcal{M}=\frac{1}{2}\left(R^2+\frac{\tilde{\Theta}^2}{r^2}\right)-\frac{\mu}{r},
\end{equation}
where
\begin{equation} \label{variedG}
\tilde{\Theta}^2=\Theta^2\left[1-\frac{1}{2}J_2\frac{\alpha^2}{p^2}(3c^2-1)\right]
\end{equation}
is function only of $\Theta$ and the parameters of the simplified problem; namely the integral $N$ and the physical parameters defining the gravitational field. The objective of the torsion is to find a canonical change of variables $(\theta,\nu,\Theta,N)\mapsto(\theta^*,\nu^*,\Theta^*,N^*)$ such that in the new variables $\tilde{\Theta}=\Theta^*$. That is, in the asterisk variables Eq.~(\ref{Hamainq}) turns into the Kepler problem Hamiltonian.
\par

Working in polar variables there is no distinction between bounded or unbounded motion, so the torsion transformation is exactly the same as in the elliptic case originally described in \S9 of \cite{Deprit1981}. Rather, we follow the notation in \cite{HautesserresLaraCeMDA2017} and denote $\Phi=\Phi(\Theta,N)\equiv\tilde\Theta/\Theta$, $\epsilon=\epsilon(\Theta)\equiv-\frac{1}{2}J_2(\alpha/p)^2<0$, and define the torsion by the sequence,
\begin{align} \label{tG}
\Theta^*= & \; \Theta\Phi, \\  \label{tz}
\theta^*= & \; \theta\Phi\left(\Phi^2-2\epsilon\frac{\partial\Phi^2}{\partial\epsilon}-\frac{1}{2}c\frac{\partial\Phi^2}{\partial{c}}\right)^{-1}, \\  \label{th}
\nu^*= & \; \nu-\frac{1}{2}\frac{\theta^*}{\Phi}\frac{\partial\Phi^2}{\partial{c}},
\end{align}
whereas not only $N$, but $r$ and $R$ remain unaltered. The inverse transformation requires first to solve $\Theta=\Theta(\Theta^*,N)$ from the  implicit function (\ref{tG}), which is then substituted in Eqs.~(\ref{tz}) and (\ref{th}) to be trivially solved in $\theta$ and $\nu$. Alternatively to the numerical solution by root-finding procedures, Eq.~(\ref{tG}) can be inverted analytically to be given in the form of a power series in $\epsilon$. If we limit to the first order of $J_2$ of the perturbation solution computed hitherto, we obtain
\begin{equation} \label{ZZfirst}
\Theta=\Theta^*\left[1-\frac{1}{2}(3c^2-1)\epsilon\right],
\end{equation}
in which, now, $c=N/\Theta^*$, and $p=\Theta^{*2}/\mu$ must be replaced in the definition of $\epsilon$. 
\par

After written in asterisk variables, Eq.~(\ref{Hamainq}) becomes a Keplerian that is next completely reduced by the exact transformation from polar to hyperbolic Delaunay variables. Thus, $\mathcal{M}=+\mu^2/(2L^2)$, from which
\[
\ell=\ell_0+nt,\quad g=g_0,\quad h=h_0,\quad L=L_0,\quad G=G_0,\quad H=H_0,
\]
with $n=-\mu^2/L^3$.
\par

Note that the Keplerian orbit obtained with the torsion may have no physical sense, showing complex values of the inclination. This fact is clearly observed in the limit case obtained when making $\Theta=N$ in Eq.~(\ref{ZZfirst}). Namely,
\[
\frac{N}{\Theta^*}\equiv{c}=1-(3c^2-1)\epsilon,
\]
a quadratic equation in $c$ that exists only in the interval $-\frac{1}{2}+\frac{1}{\sqrt{6}}\le\epsilon<0$, in which $c>1$ always. Irrespective of the meaning that $c$ may have in Eq.~(\ref{ZZfirst}) in these particular cases, the torsion transformation still remains valid.
\par

\section{Higher-order extensions}

In the case of bounded orbits, a second order solution of the main problem based on a quasi-Keplerian Hamiltonian would be unavoidably harmed by the appearance of secular terms in the transformation equations. However, these kinds of terms are not at all of worry in the case of unbounded motion, in which their range of variation is constrained by the limit values of the asymptotes of a hyperbola. Therefore, we can safely extend the computed first order solution to higher orders following exactly the same approach. Namely, first computing the mapping that yields a quasi-Keplerian Hamiltonian while accepting the boundary conditions at infinity, and then applying a torsion to obtain the final Keplerian.
\par

Thus, from the known terms obtained at the second order of the elimination of the parallax, we select those terms that are free from the true anomaly but also from the argument of the periapsis and the eccentricity. The latter, we recall, is a function of the radius and radial velocity in addition to the specific angular momentum. We obtain
\begin{equation} \label{M02}
\mathcal{M}_{0,2}=-\frac{\Theta^2}{r^2}\frac{\alpha^4}{p^4}\frac{1}{16}(21s^4-42s^2+20),
\end{equation}
The simple consideration of this additional term in the reduction of the quasi-Keplerian Hamiltonian $\mathcal{M}=\mathcal{M}_{0,0}+J_2\mathcal{M}_{0,1}+\frac{1}{2}J_2^2\mathcal{M}_{0,2}$ yields refinements in the computation of the final Keplerian mean motion, as well as in the solution of the torsion transformation given by Eqs.~(\ref{tG})--(\ref{th}), where, now,
\begin{equation} \label{Phi2nd}
\Phi=\left[1-\epsilon(3c^2-1)+\frac{1}{4}\epsilon^2(1-21c^4)\right]^{1/2}.
\end{equation}
In addition, the analytical inversion of Eq.~(\ref{tG}) provided by Eq.~(\ref{ZZfirst}) must be replaced by the second order of $J_2$ solution given by
\begin{equation} \label{G2nd}
\Theta=\Theta^*\left[1-\frac{\epsilon}{2}(3c^2-1)-\frac{3\epsilon^2}{4}(2c^2-1)(5c^2-1)\right],
\end{equation}
where $c$, and $p$ are now functions of the asterisk variables.
\par

Major improvements will be obtained if, besides, we complement the analytical solution with second order corrections to the elimination of the parallax. That is, Eq.~(\ref{tx1}) is replaced by
\begin{equation} \label{tx2}
\xi=\xi'+J_2\xi_{0,1}|_{\xi=\xi'}+\frac{1}{2}J_2^2\xi_{0,2}|_{\xi=\xi'},
\end{equation}
where $\xi_{0,2}\equiv\{\xi_{0,1},\mathcal{U}_1\}+\{\xi,\mathcal{U}_2\}$. To do that, we need to compute the generating function term $\mathcal{U}_2$ of the Lie transformation that yields $\mathcal{M}_{0,2}$. It is obtained in the usual way up to an arbitrary constant $Q_0\equiv{Q}_0(g,h,L,G,H)$. 
\par

The new corrections need to match the boundary condition at infinity too. Analogously to the first order, we readily find that $R_{0,2}|_{\infty}\equiv0$, and hence the boundary conditions are automatically fulfilled by the corrections to the radial velocity. On the contrary, the equation $\Theta_{0,2}|_{\infty}=0$ involves the partial derivative $\partial{Q}_0/\partial{g}$, which is solved by indefinite integration to give
\begin{align} \nonumber
Q_0= & \; G\frac{\alpha^4}{p^4}\frac{1}{256 e^2}\Big\{
-2s^2\left[ 3e^4(17s^2-18)+8e^2(75 s^2-68)+8(s^2+6) \right. \\ \nonumber
& \left. +96\eta^3(5s^2-4)(\pi+\arctan\eta) \right]\sin2g
-3(3e^4+6e^2-16)s^4\sin4g \\ \nonumber
& +\frac{4s^2}{\eta^2}\left\{ 6\eta^2\left[e^4(15s^2-14)+4(3 e^2-2)(5s^2-4)\right](\pi+\arctan\eta)
\right. \\ \nonumber
& \left. +\eta\left[e^4(278-329s^2)+e^2(298s^2-284)+4(s^2+6)\right]\right\}\cos2g \\ \nonumber
& +6\eta(e^2+8)s^4\cos4g \Big\}+Q_1(L,G,H).
\end{align}
The new integration constant $Q_1$ is solved from the equation $\nu_{0,2}|_{\infty}=0$, which only involves the partial derivative $\partial{Q}_1/\partial{H}$. We obtain
\begin{align} \nonumber
Q_1=&\; G\frac{\alpha^4}{p^4}\frac{3}{128\eta}(1-s^2)\left[2e^2\eta(5s^2+13)(\pi+\arctan\eta)\right.\\ \nonumber
& \left. -e^2(11s^2+75)-85s^2+107\right] +Q_2(L,G).
\end{align}
Analogously, the integration constant $Q_2$is solved from the equation $\theta_{0,2}|_{\infty}=0$, which only involves the partial derivative $\partial{Q}_2/\partial{G}$. We obtain
\begin{align} \nonumber
Q_2= &  \; G\frac{\alpha^4}{p^4}\frac{1}{128\eta}\left[81e^2-30e^2\eta(\pi+\arctan\eta)-49\right]+Q_3(L).
\end{align}
We take $Q_3\equiv0$ for simplicity, and finally obtain
\begin{align} \nonumber
\mathcal{U}_2= & \; G\frac{\alpha^4}{p^4}\frac{3}{64e^2}\left\{
\left[2e^4(15s^2-14)+8(3e^2-2)(5s^2-4)\right]s^2\cos2g \right. \\ \nonumber
& \left. -16\eta^3(5s^2-4)s^2\sin2g -e^4(5s^4+8s^2-8) \right\}\psi
\\ \nonumber
& 
+G\frac{\alpha^4}{p^4}\frac{1}{256e^3\eta}\sum_{k=0}^2s^{2k}
\bigg\{\sum_{j=j_0}^5\sum_{i=0}^3q_{k,i,j}e^{2i+1-(j\bmod2)}\cos(jf+2kg) \\ \label{U2}
& +\eta\sum_{j=j_0}^{j_1}\sum_{i=0}^2p_{k,i,j}e^{2i+1-(j\bmod2)}\sin(jf+2kg)
\bigg\},
\end{align}
where $j_0=-2 (k \bmod 2)-1$, $j_1=6-2 (k \bmod 2)$, $\psi=\pi-f+\arctan\eta$, and the inclination polynomials $p_{k,i,j}$ and $q_{k,i,j}$ are listed in Tables \ref{t:incpolcos} and \ref{t:incpolsin}. The second order corrections $\xi_{0,2}$ to the polar variables are then obtained from the straightforward evaluation of the corresponding Poisson brackets. 
\par

Note that the appearance of the function $\psi$, which also exists in Hori's first order solution, was only delayed to the second order of our approach. This postponement resembles analogous effects in the case of bounded motion, in which the elimination of the parallax simplification delays by one order of the perturbation approach the need of integration by parts of terms involving the equation of the center in the computation of the generating function \cite{Healy2000,Lara2020arxiv}.
\par

\begin{table}[htbp]
\centering \small
\begin{tabular}{@{}rccc@{}}
$i,j$ & $k=0$ & $k=1$ & $k=2$ \\
\hline\noalign{\smallskip}
$_{0,-2}$ & $0$ & $64(3 s^2-2)$ & $0$ \\
$_{0,-1}$ & $0$ & $160 (3 s^2-2)$ & $0$ \\
$_{0,0}$ & $0$ & $16 (s^2+6)$ & $-48$ \\
$_{0,1}$ & $0$ & $0$ & $-48$ \\
$_{0,2}$ & $-24 s^2 (15 s^2-4)$ & $0$ & $-120$ \\
$_{0,3}$ & $-112 s^4$ & $0$ & $0$ \\
$_{0,4}$ & $-72 s^4$ & $0$ & $0$ \\
$_{1,-3}$ & $0$ & $8 (3 s^2-2)$ & $0$ \\
$_{1,-2}$ & $0$ & $-80 (3 s^2-2)$ & $0$ \\
$_{1,-1}$ & $-12 s^2 (13 s^2-4)$ & $-8 (45 s^2-26)$ & $-12$ \\
$_{1,0}$ & $2 (255 s^4-576 s^2+272)$ & $8 (149 s^2-142)$ & $42$ \\
$_{1,1}$ & $-12 s^2 (13 s^2-4)$ & $96 (2 s^2-1)$ & $-36$ \\
$_{1,2}$ & $12 s^2 (45 s^2-16)$ & $-96 (6 s^2-5)$ & $60$ \\
$_{1,3}$ & $4 s^2 (7 s^2+8)$ & $0$ & $-120$ \\
$_{1,4}$ & $90 s^4$ & $0$ & $-54$ \\
$_{1,5}$ & $-12 s^4$ & $0$ & $0$ \\
$_{2,-3}$ & $0$  & $-10 (3 s^2-2)$ & $0$ \\
$_{2,-2}$ & $0$ & $16 (3 s^2-2)$ & $0$ \\
$_{2,-1}$ & $8 (75 s^4-72 s^2+20)$ & $-2 (153 s^2-134)$ & $15$ \\
$_{2,0}$ & $6 (11 s^4+64 s^2-48)$ & $-4 (329 s^2-278)$ & $6$ \\
$_{2,1}$ & $8 (75 s^4-72 s^2+20)$ & $-6 (127 s^2-90)$ & $81$ \\
$_{2,2}$ & $4 (27 s^4-72 s^2+32)$ & $24 (9 s^2-10)$ & $60$ \\
$_{2,3}$ & $s^2 (145 s^2-64)$ & $-2 (195 s^2-146)$ & $105$ \\
$_{2,4}$ & $-18 s^4$ & $-36 (3 s^2-2)$ & $54$ \\
$_{2,5}$ & $15 s^4$ & $0$ & $-9$ \\
$_{3,-3}$ & $0$ & $2 (3 s^2-2)$ & $0$ \\
$_{3,-1}$ & $-6 (5 s^4+4 s^2-4)$ & $24 (7 s^2-6)$ & $-3$ \\
$_{3,1}$ & $-6 (5 s^4+4 s^2-4)$ & $12 (25 s^2-22)$ & $3$ \\
$_{3,3}$ & $-25 s^4-16 s^2+16$ & $8 (15 s^2-14)$ & $15$ \\
$_{3,5}$ & $-3 s^4$ & $-6 (3 s^2-2)$ & $9$ \\
\noalign{\smallskip}\hline
\end{tabular}
\caption{Inclination polynomials $q_{k,i,j}$ in Eq.~(\protect\ref{U2}).}
\label{t:incpolcos}
\end{table}
\begin{table}[htb]
\centering \small
\begin{tabular}{@{}rccc@{}}
$i,j$ & $k=0$ & $k=1$ & $k=2$ \\
\hline\noalign{\smallskip}
$_{0,-2}$ & $0$ & $ -64 (3 s^2-2)$ & $0$ \\
$_{0,-1}$ & $0$ & $ -160 (3 s^2-2)$ & $0$ \\
$_{0,0}$ & $0$ & $ -16 (s^2+6)$ & $48$ \\
$_{0,1}$ & $0$ & $0$ & $48$ \\
$_{0,2}$ & $-24 s^2 (15 s^2-4)$ & $0$ & $120$ \\
$_{0,3}$ & $-112 s^4$ & $0$ & $0$ \\
$_{0,4}$ & $-72 s^4$ & $0$ & $0$ \\
$_{1,-3}$ & $0$ & $ -8 (3 s^2-2)$ & $0$ \\
$_{1,-2}$ & $0$ & $ 48 (3 s^2-2)$ & $0$ \\
$_{1,-1}$ & $12 s^2 (13 s^2-4)$ & $ 24 (5 s^2-2)$ & $12$ \\
$_{1,0}$ & $0$ & $ -16 (75 s^2-68)$ & $-18$ \\
$_{1,1}$ & $-12 s^2 (13 s^2-4)$ & $ -96 (2 s^2-1)$ & $60$ \\
$_{1,2}$ & $72 s^2 (5 s^2-2)$ & $ -32 (15 s^2-13)$ & $0$ \\
$_{1,3}$ & $4 s^2 (8-7 s^2)$ & $0$ & $120$ \\
$_{1,4}$ & $54 s^4$ & $0$ & $66$ \\
$_{1,5}$ & $-12 s^4$ & $0$ & $0$ \\
$_{2,-3}$ & $0$ & $ 6 (3 s^2-2)$ & $0$ \\
$_{2,-1}$ & $-2 (27 s^4+180 s^2-128)$ & $ 18 (17 s^2-14)$ & $-9$ \\
$_{2,0}$ & $0$ & $ -6 (17 s^2-18)$ & $-9$ \\
$_{2,1}$ & $2 (27 s^4+180 s^2-128)$ & $ -2 (825 s^2-742)$ & $-45$ \\
$_{2,2}$ & $6 (5 s^4+8 s^2-8)$ & $ -24 (s^2-2)$ & $-15$ \\
$_{2,3}$ & $3 s^2 (39 s^2-16)$ & $ 2 (55 s^2-34)$ & $-15$ \\
$_{2,4}$ & $0$ & $ 6 (13 s^2-10)$ & $-3$ \\
$_{2,5}$ & $9 s^4$ & $0$ & $21$ \\
$_{2,6}$ & $0$ & $0$ & $3$ \\
\noalign{\smallskip}\hline
\end{tabular}
\caption{Inclination polynomials $p_{k,i,j}$ in Eq.~(\protect\ref{U2}).}
\label{t:incpolsin}
\end{table}

\section{Performance of the analytical solution}
%{\color{red} POR SI LO PREGUNTA ALGUN REVISOR
%Practical use of the perturbation theory for astrodynamics applications is under investigation and will be approached in a forthcoming paper. Here we only provide three examples that serve to illustrate the performance that can be expected from the new analytical solution in different scenarios.
%}
%\par

For reference, in the first two examples we choose the physical parameters of the Earth ($\mu=398600.44$ $\mathrm{km^3/s^2}$, $\alpha=6378.1363$ km, $J_2=0.001082634$). The altitude of Earth flybys for gravity assist typically range from about 300 km, as in the case of NASA's Galileo mission,\footnote{\href{https://solarsystem.nasa.gov/missions/galileo/in-depth/}{solarsystem.nasa.gov/missions/galileo/in-depth} Retrieved July 18, 2022} to the more common several thousands of km, as for instance for NASA's Stardust mission.\footnote{\href{https://solarsystem.nasa.gov/stardust/news/ega/index.html}{solarsystem.nasa.gov/stardust/news/ega} Retrieved July 18, 2022} % https://www.wikiwand.com/en/List_of_Earth_flybys
Accordingly, we fix a perigee height of 1000 km over the Earth's surface in the examples, which is low enough to clearly undergo the oblateness disturbing effect.
\par

The first test case is for a high eccentricity hyperbolic-type orbit, for which the flyby happens in a relatively short time, and both the Kepler approach and the first order of the perturbation solution are expected to provide acceptable results. The second example deals with a quasi-parabolic orbit, for which the spacecraft remains much more time in the close vicinity of the Earth, and, therefore, it is expected to be a more challenging test for the analytical solution.
\par

\begin{figure}[htb]
\centerline{
\includegraphics[scale=0.8]{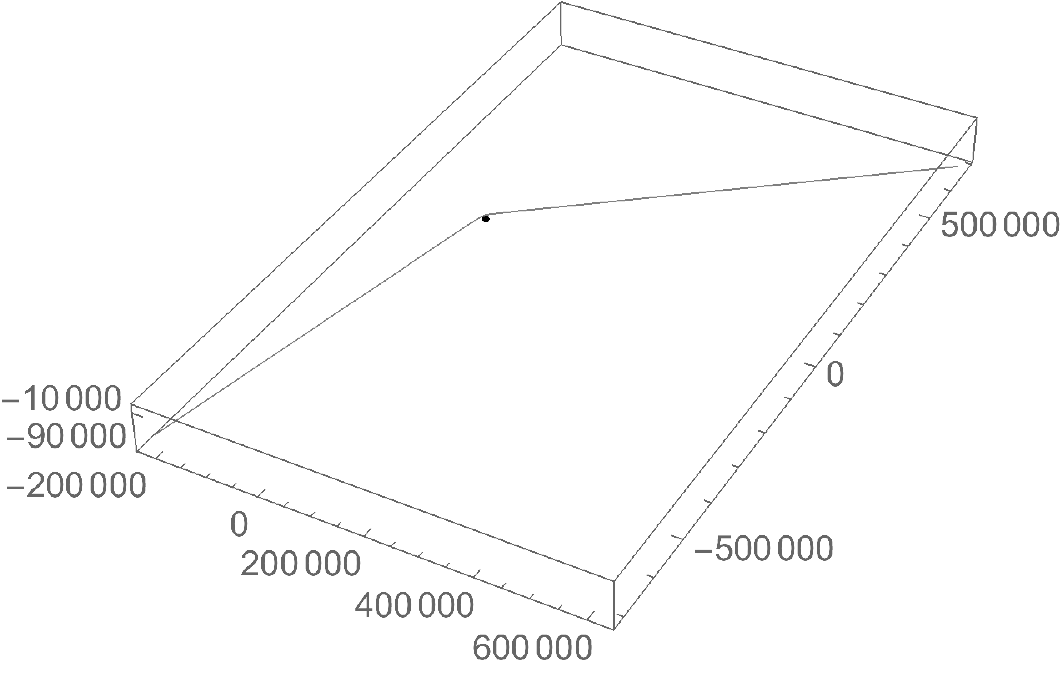}
\includegraphics[scale=0.7]{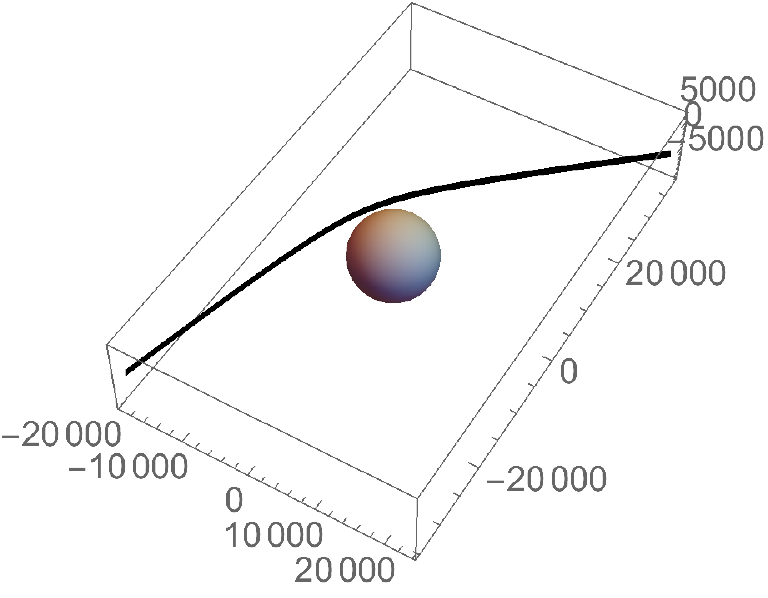}
}
\caption{Earth's flyby with $e=4$ and 1000 km periapsis distance. Left: true orbit. Right: magnification on the closest approach. Distance in km.}
\label{f:fb1}
\end{figure}

In the first case we choose initial conditions corresponding to the orbital elements $a=2459.38$ km, $e=4$, $I=23.5^\circ$, $\Omega=60^\circ$, $\omega=90^\circ$, $M=-21400^\circ$, that correspond approximately to the entrance of the satellite into the Earth's sphere of influence (SOI). The ``true'' reference solution is computed from a numerical integration of these initial conditions in the main problem dynamics that guarantees the preservation of both the energy integral and the third component of the angular momentum within at least 14 digits along the propagation interval of 36 hours. The reference orbit is shown in Fig.~\ref{f:fb1}.
\par

\begin{figure}[htb]
\centerline{
\includegraphics[scale=0.8]{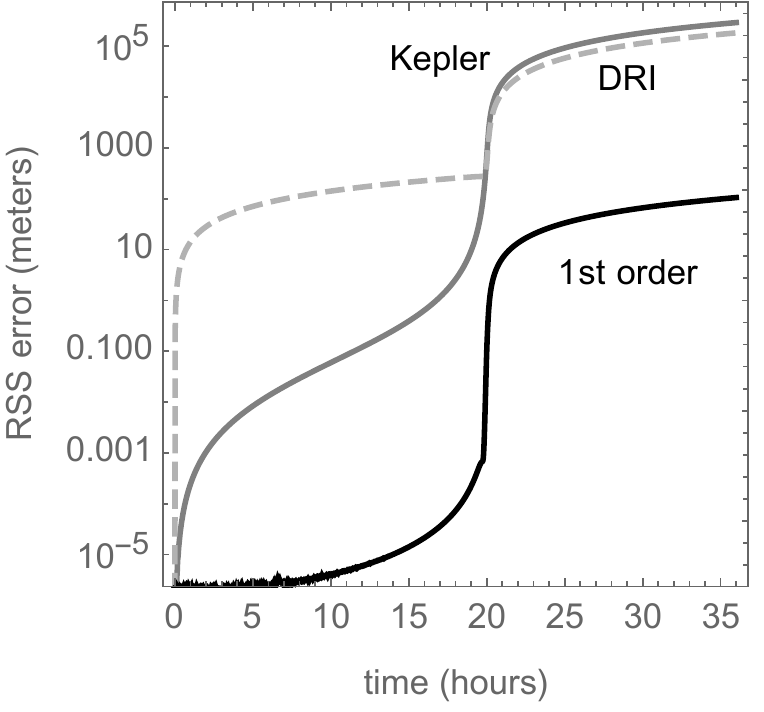}
\includegraphics[scale=0.8]{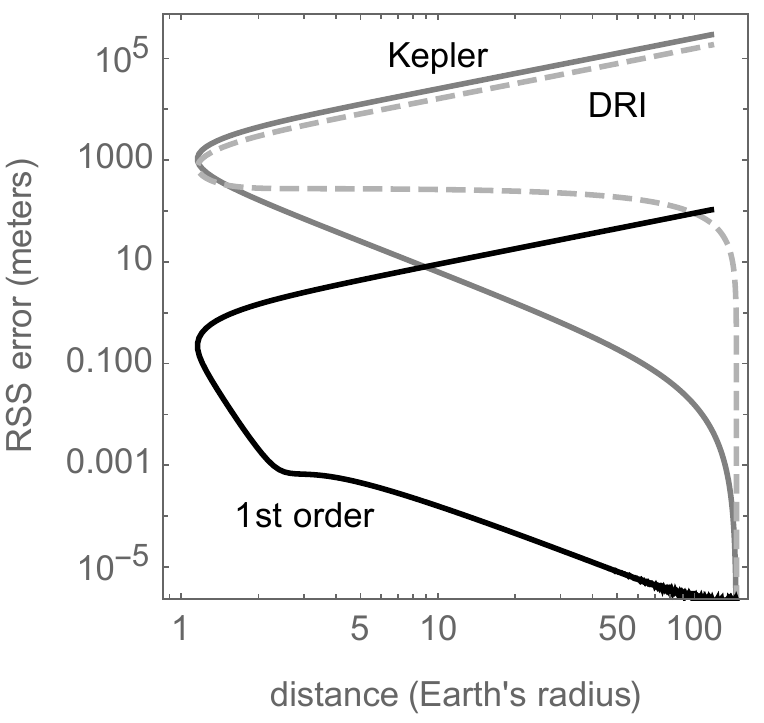}
}
\caption{RSS position errors of the flyby in Fig.~\protect\ref{f:fb1} when using the Keplerian dynamics, DRI-common, and the 1st order of the perturbation solution. Note the logarithmic scale.}
\label{f:fb1err}
\end{figure}

Then, the initial conditions are propagated analytically, first in the Keplerian dynamics, next, in the ``common'' interpretation of Deprit's radial intermediary, in what follows abbreviated DRI, in which Eqs.~(\ref{Hamainq})--(\ref{variedG}) are assumed to remain valid in osculating elements, and finally in the ``natural'' interpretation of DRI \cite{Deprit1981}; that is, the approximation of the main problem dynamics provided by the new perturbation solution. The Root Sum Square (RSS) errors of each analytical solution with respect to the true orbit are shown in Fig.~\ref{f:fb1err}. The Keplerian propagation starts with zero error, and, because the disturbing effect of $J_2$ is very small at SOI entrance, errors grow slowly, reaching the meter level at a distance of approximately 25 Earth's radius, just a few hours before the closest approach. Errors grow notably during the perigee passage, and, due to the different mean Keplerian motion and the one of the main problem dynamics, increase almost linearly after that. At the end of the propagation, the RSS errors of the Keplerian approximation reach several hundreds of km.
\par

When using the ``common'' version of DRI, we get closer to the $J_2$ dynamics although of a slightly different orbit due to the lack of transformation between original and prime variables. By this reason, the errors in the arrival branch of the predicted orbit grow faster than the Keplerian approximation until the accumulation of nonlinear effects due to the $J_2$ perturbation is significant enough, which happens at about 1.5 Earth's radius for this orbit. Then, the better representation of the main problem dynamics by DRI results in clearly smaller errors than those of the Keplerian approximation during the flyby, and also in the departure branch of the orbit, as illustrated with the dashed gray lines in both plots of Fig.~\ref{f:fb1err}. The computational burden of DRI is notably alleviated when the root finding procedure needed in the solution of the implicit function Eq.~(\ref{tG}) is replaced by the direct evaluation of Eq.~(\ref{ZZfirst}), yet in this case the analytical propagation starts with an initial error of half a meter ---which is clearly irrelevant at the precision of the analytical propagation.
\par

Things notably improve when using the natural version of DRI, labelled ``1st order'' in Fig.~\ref{f:fb1err}. Now we unavoidably deal with the mandatory transformation of the initial osculating state to the mean elements that initialize the constants of the perturbation solution ---using the corrections in Eqs.~(\ref{r01})--(\ref{G01}) with opposite signs--- and the consequent recovery of osculating elements using Eqs.~(\ref{r01})--(\ref{G01}) in order to compute the RSS error.\footnote{Recall that, while direct and inverse transformations are just opposite for first order corrections, they are fed with different sets of elements, either mean or osculating, and hence their composition is affected by the truncation error of the perturbation theory} Errors provided by the 1st order solution grow fast at perigee passage, but their following almost linear increase is, approximately, $J_2$ times smaller than in the Keplerian and DRI cases due to the much better approximation provided by the perturbation solution, and barely reaches 100 m at the end of the propagation.
\par

\begin{figure}[htb]
\centerline{
\includegraphics[scale=0.8]{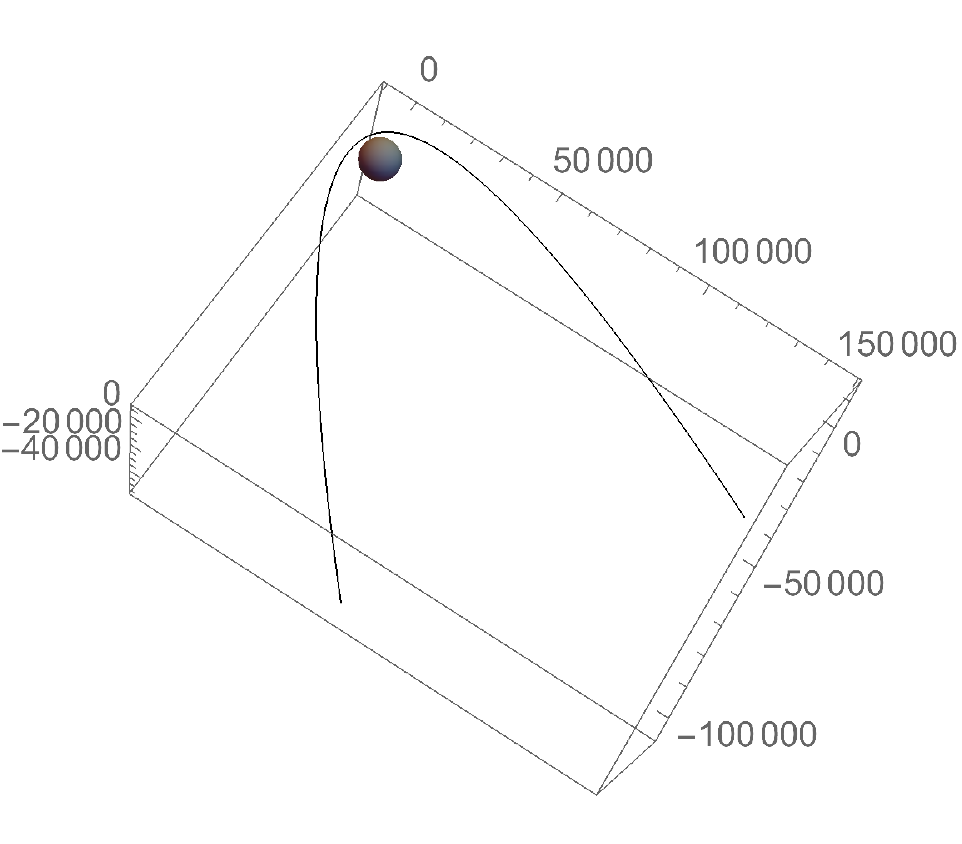}
}
\caption{Quasi-parabolic flyby 1000 km above Earth's surface. Distance in km.}
\label{f:fb2}
\end{figure}

For the second case we choose the quasi-parabolic orbit $a=1.47563\times 10^6$ km, $e=1.005$, $I=23.5^\circ$, $\Omega=60^\circ$, $\omega=90^\circ$, $M=-1^\circ$, which will fly about 1000 km over the Earth's surface too, but now during a much longer time than in the previous example. The reference orbit, which is shown in Fig.~\ref{f:fb2}, was numerically integrated in the main problem dynamics up to 24 hours. As before, it is compared with equivalent solutions provided by the Keplerian approximation, on the one hand, and the torsion-based solution of the main problem, both in the common and natural versions, on the other. 
\par

\begin{figure}%[htb]
\centerline{
\includegraphics[scale=0.8]{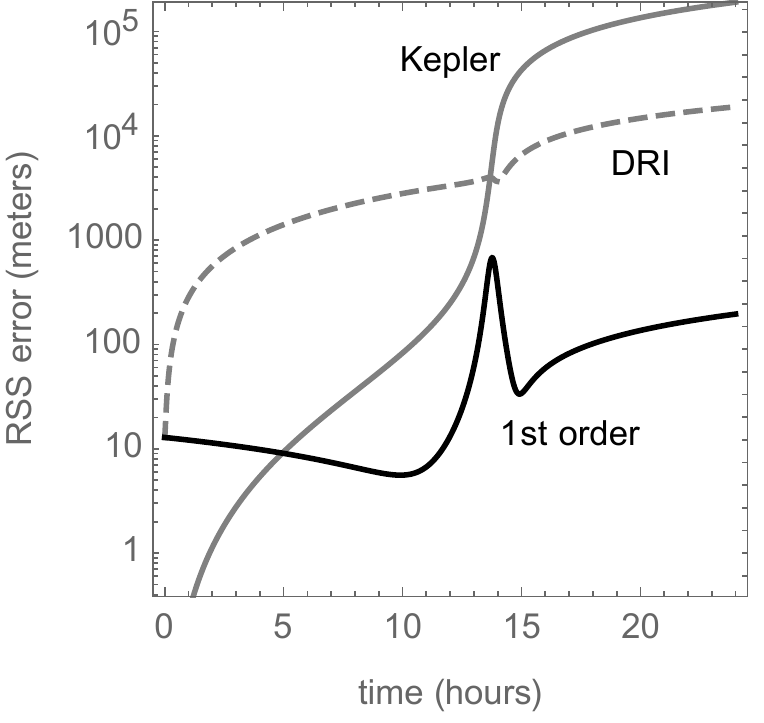}
\includegraphics[scale=0.8]{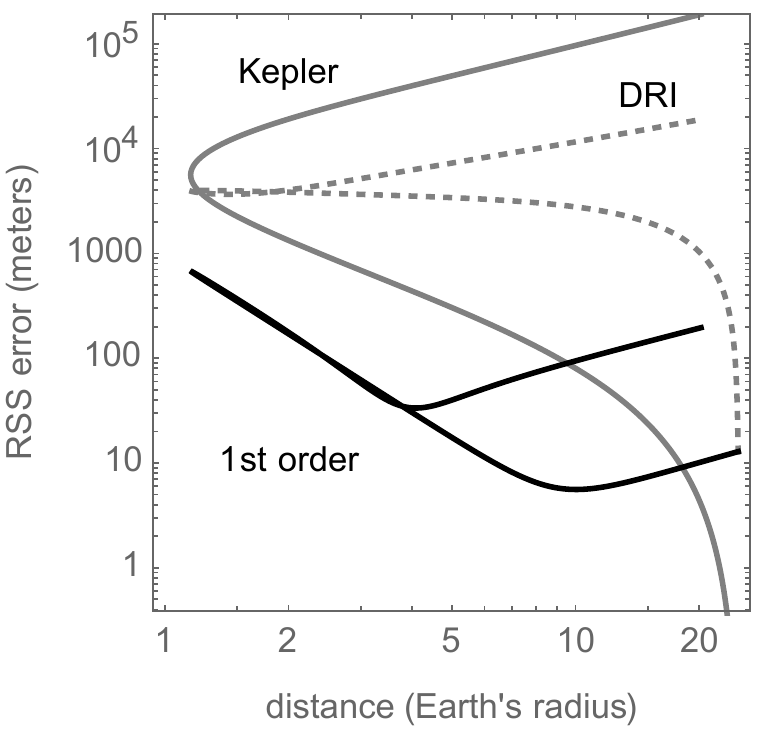}
}
\caption{RSS position errors of the flyby in Fig.~\protect\ref{f:fb2} when using the Keplerian dynamics and the 1st order of the perturbation solution. Note the logarithmic scale.}
\label{f:fb2err}
\end{figure}

Now the spacecraft undergoes the disturbing effects of the oblateness perturbation for a much longer time than in the previous flyby, and hence the improvements of DRI over the Keplerian approximation are much more evident in this example, as shown in Fig.~\ref{f:fb2err} ---where the approximate inversion in Eq.~(\ref{ZZfirst}) has been used and hence the apparent initial error of the torsion-based solution. Also, analogous improvements to the previous example are also obtained by the 1st order perturbation solution in this case, for which the errors remain $\mathcal{O}(J_2)$ with respect to corresponding ones of the Kepplerian approximation. However, as shown in Fig.~\ref{f:fb2err}, the errors of the first order solution now clearly peak along the close visit to Earth reaching a RSS error of about 700 m at perigee (about 10\% of the minimum distance to the Earth's surface), which immediately falls by one order of magnitude, finally reaching about 200 m at the end of the propagation. The latter figure means a $\mathcal{O}(J_2)$ improvement when compared to about the same number of km, instead of m, provided by the Keplerian approach. The peak, we recall, is due to the appearance of $\eta$ in denominators of some of the first order corrections, which in this more challenging case takes the value $\eta\approx0.1$. On the other hand, it deserves mentioning that the fact that the errors of the 1st order solution decrease form the initial 10 meter level to just a few meters when approaching Earth is just an apparent paradox that is easily explained by the approximate truncation of the inversion of the torsion to the first order of $J_2$, and the decreasing value of $\theta$ in Eq.~(\ref{tz}). Eventually, errors in the dynamical modeling provided by the natural version of DRI overcome those due to the approximation in Eq.~(\ref{ZZfirst}), and grow again irrespective of the $\theta$ value.
\par

\begin{figure}%[htb]
\centerline{
\includegraphics[scale=0.9]{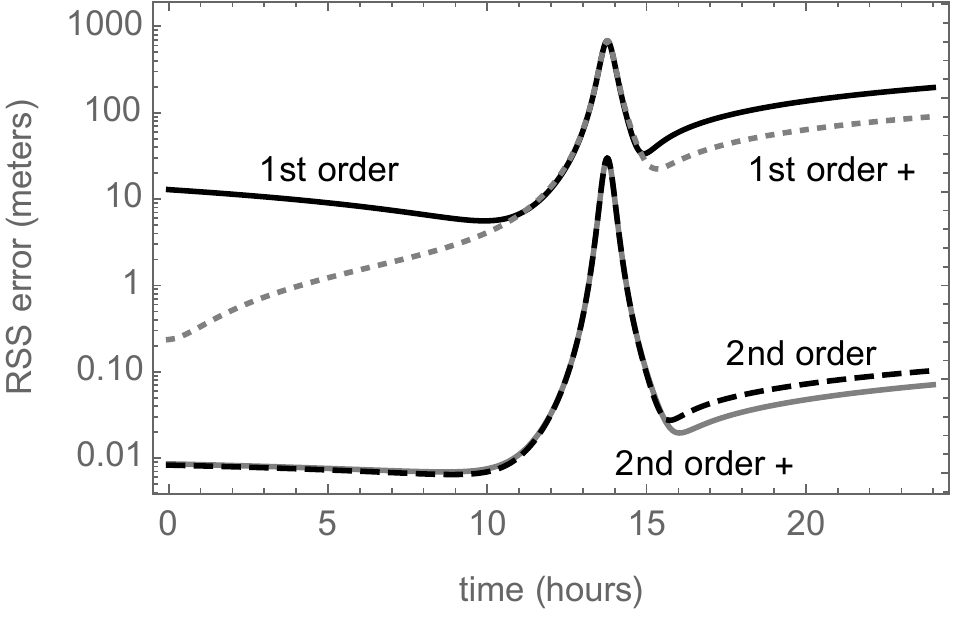}
}
\caption{RSS position errors of the flyby in Fig.~\protect\ref{f:fb2} when using higher orders of the perturbation solution. Note the logarithmic scale.}
\label{f:fb2er2}
\end{figure}

Like in the case of bounded motion, a general improvement of the 1st order solution is obtained when including second-order terms in the mean elements propagation, a strategy that negligibly increases the computational burden \cite{Brouwer1959,CoffeyAlfriend1984}. In our case, this means including the term in Eq.~(\ref{M02}) in the quasi-Keplerian Hamiltonian, as well as Eqs.~(\ref{Phi2nd}) and (\ref{G2nd}) in the solution of the torsion transformation. As shown in Fig.~\ref{f:fb2er2}, while these additional terms generally reduce the errors of the analytical propagation (labelled ``1st order +''), they do not show improvements during the close approach to the Earth. Certainly, computing the second order corrections derived from Eq.~(\ref{U2}) makes a full second order solution that, as expected, generally improves the accuracy to $\mathcal{O}(J_2^2)$, as is also shown in Fig.~\ref{f:fb2er2} (black-dashed curve labelled ``2nd order''). The fact that the expected improvement does not happen at perigee to the whole extent, is due to the appearance in denominators of the second order corrections of higher powers of $\eta$ than in the first order approach, which now reach the fourth power in the case of $\theta_{0,2}$ and $R_{0,2}$. Still, the errors of the 2nd order solution at perigee passage remain at least 20 times smaller than those of the first order solution. Finally, the RSS error curve labelled ``2nd order +'' includes third-order secular terms obtained from the elimination of the parallax in the perturbation solution. Like in the first order case, these additional terms slightly improve the computation of the ``mean'' mean motion ---the mean motion of the quasi-Keplerian Hamiltonian--- yet their effects are only noted in the departure branch of the orbit.
\par

Because perturbation theories of this kind are rather uncommon, the final example is intended to show that the convergence of the successive orders of the perturbation solution, which is commonly observed in the case of bounded perturbed Keplerian motion (see \cite{LaraPalacianRussell2010,Lara2020}, for instance), also happens in the current case of perturbed hyperbolic Keplerian motion. To do that, we test our solution with the challenging case of a fish-type orbit \cite{MartinusiGurfil2013}. However, since parabolic orbits are necessarily excluded from the current perturbation solution due to the $\eta$ divisor that affects the elimination of the parallax transformation, we rather base our tests on an equivalent hyperbolic branch of the Newtonian rosetta.
\par

\begin{figure}[htb]
\centerline{
\includegraphics[scale=0.7]{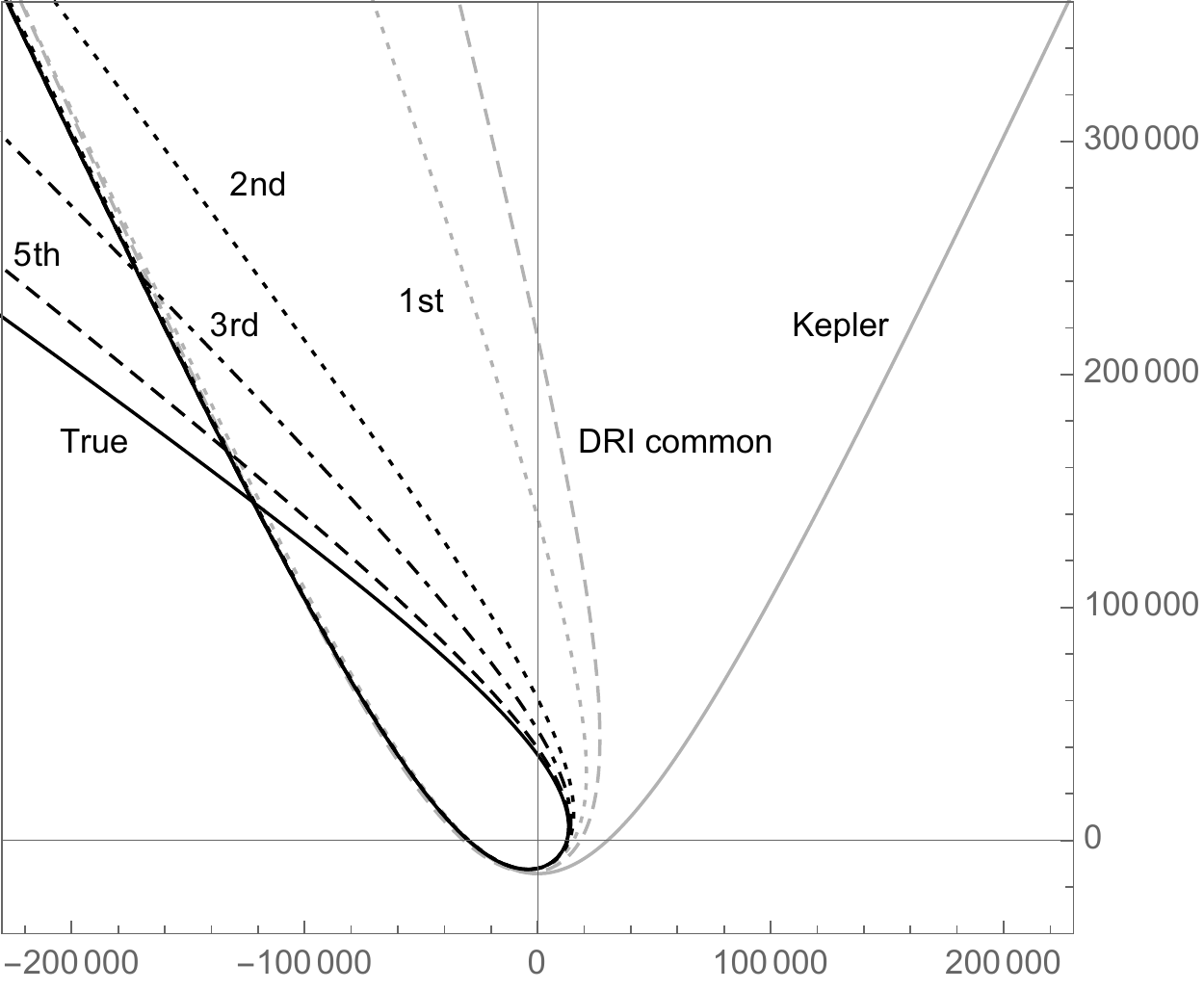}
}
\caption{True hyperbolic rosetta branch (black curve) and successive approximations of the perturbation solution truncated to the zeroth (Keplerian, gray), first (gray, dotted), second (black, dotted), third (black, dot-dashed) and fifth order (black, dashed) of Jupiter's main problem. The gray-dashed curve was generated with the ``common'' version of DRI. Distances are km.}
\label{f:fish}
\end{figure}

Thus, for the gravitational parameters of Jupiter ($\mu=1.268\times10^8\,\mathrm{km^3/s^2}$, $\alpha=71492$ km, $J_2=0.01475$) we choose the initial conditions of an equatorial orbit with $a=724920$ km, $e=1.1$, $g=270^\circ$, and $\ell=-120^\circ$. The fact that this orbit would impact Jupiter does not worry us in our aim of validating the perturbation theory. Results are summarized in Fig.~\ref{f:fish}, where the true orbit, computed numerically, is represented with a black curve, approaching to Jupiter from the top-left corner, surrounding the attraction center by more than $180^\circ$ to first cross the approaching branch of the orbit and then leave the scene by the left side of the picture. The Keplerian approximation (gray curve of Fig.~\ref{f:fish}) only mimics the true trajectory in the approaching part, and then takes a clearly divergent path. DRI (gray dashed curve) is able to turn the departure branch towards the true solution, a result that is clearly improved by the first order theory (gray dotted curve), yet not enough enhanced to produce the crossing between the approaching and departing branches of the orbit. This intersection is obtained with the second order solution (black dotted curve of Fig.~\ref{f:fish}), yet it happens out of the frame of Fig.~\ref{f:fish}, clearly far away from the true event. The third-order approximation (black dot-dashed curve of Fig.~\ref{f:fish}) notably improves the results, and the crossing now occurs inside the displayed area in Fig.~\ref{f:fish}. Finally, the black-dashed curve, which was computed with a fifth-order truncation of the elimination of the parallax complemented with additional secular terms coming from a sixth-order truncation of the quasi-Keplerian Hamiltonian, reasonably approximates the true orbit. It is expected that higher order truncations continue to converge slowly to the true solution, yet we did not compute them due to the important growth in the size of the perturbation series and the corresponding exponential increase of the computational burden. We are satisfied with these results, which show the correctness of Hori's approach for dealing with unbounded perturbed Keplerian motion.
\par

%\section{Hyperbolic fish-type orbits}

\section{Remarks on the boundary condition: the third body perturbation case}

Similarly to Hori's approach, the transformation that yields Deprit's radial intermediary depends on the choice made for the integration constant. The main problem straightforwardly suggests to impose null perturbing effects at infinity, as at great distances from the primary the non-uniform gravity field of the planet flown by is negligible.
\par

Yet, in practical examples for the Solar or a planetary system, the great distances from the flyby body are characterized by stronger effects of the main attractor of the whole system. We therefore proceed at assessing whether perturbation approaches of this kind, with Keplerian boundary condition, are suitable to study the perturbed hyperbolic Kepler problem with third body effects only. 
\par

Identifying with the subscript $B$ the perturbing body, the full Hamiltonian of the third body perturbation case is
\begin{equation}\label{eq:ThirdBodyHamFull}
	\mathcal{M} = \mathcal{M}_0 + \mathcal{M}_B = -\frac{\mu}{2a} - \mu_B \bigg(\frac{1}{|\mathbf{r}-\mathbf{r}_B|}-\frac{\mathbf{r}\cdot\mathbf{r}_B}{r_B^3}\bigg)
\end{equation}
\par

Following Lidov and Kozai \cite{Lidov1961,Kozai1962}, we replace $\mathbf{r}\cdot\mathbf{r}_B = r r_B\cos\psi$ and model the perturbing potential $\mathcal{M}_B$ as a series expansion on the Legendre polynomials of $\cos\psi$:
\begin{equation}\label{eq:ThirdBodyHamExp}
	\mathcal{M}_B = -\frac{\mu_B}{r_B} \sum_{j=2}^{+\infty} \frac{\varepsilon^j}{j!} j! \, \frac{r^j}{r_B^j} P_j(\cos\psi)
\end{equation}
For instance, $P_0(x)=1$, $P_1(x)=x$, $P_2(x)=\frac{1}{2}(3x^2-1)$, $P_3(x)=\frac{1}{2}(5x^3-3x)$, and so on. The small parameter $\varepsilon$ is just a token to identify the strength of the perturbation effects with the ratio $r/r_B$, that is, retaining few orders only provide an accurate approximation if $r\ll{r}_B$. More generally, the series converges if $r<r_B$. Note the absence of the first two orders, consequence of $P_0$ leading to the constant expression $\mu_B/r_B$ and the first order term canceling with the tidal part of the perturbing potential.
\par

Following the same principle of Hori's approach for the main problem, we seek a condition where the third body perturbation vanishes by itself. We approximate
\begin{equation}\label{eq:thirdbodyapprox}
	\mathcal{M}_B \approx -\frac{\mu_B}{r_B} 2 \, \frac{r^2}{r_B^2} \frac{1}{2}\big(3\cos^2\psi -1\big)
\end{equation}
and observe that $\cos^2\psi=1/3$ makes the second term vanish. We could be tempted to use this particular point to set the Keplerian boundary condition, however, this unique setting corresponds to the second order term only being null, rather than the full original potential. In addition, trajectories such that $\cos^2\psi=1/3$ is never verified may exist. Considering the higher order terms, we do not find a value of $\cos\psi$ that makes all the terms in the expansion to vanish.
\par

\begin{figure}[htb]
\centerline{
\includegraphics[width=\textwidth]{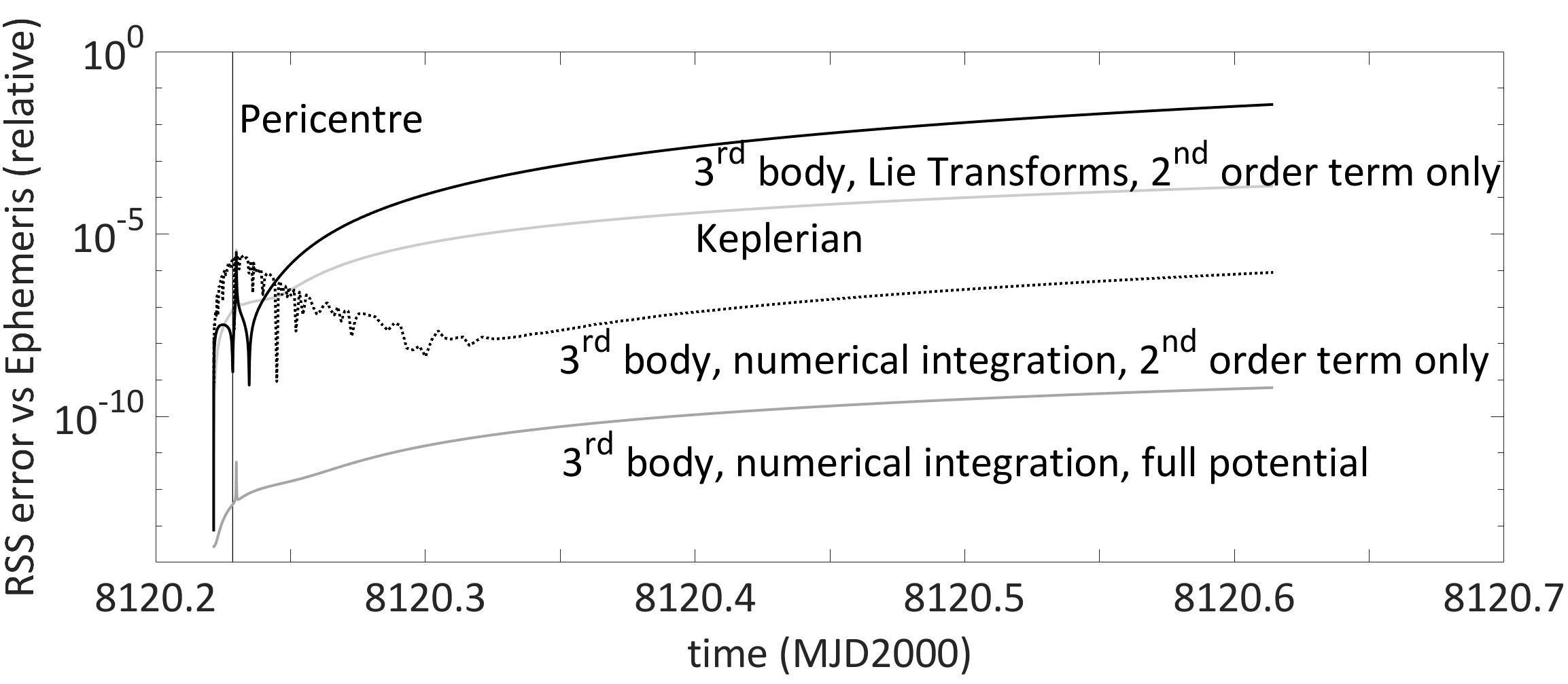}
}
\caption{Evolution of the error for the third body models, relative to the ephemeris solution.}
\label{f:ThirdBodyError}
\end{figure}

These considerations suggest that the boundary condition derived from a Keplerian hyperbola may not be suitable for the third body perturbation case. This aspect remarks that this method remains valid for the main problem only, or more in general for vanishing perturbations at infinity. We strengthen our claim with a numerical proof of the just mentioned observation, whose analytical equations of motion in prime variables we obtained with the help of a symbolic processor, following the same steps of Hori's approach but applying the Lie transform method. We do not include them in this article for the sake of brevity. We study the case of a a planetary flyby of Venus modeling the third body perturbation of the Sun, within Venus' SOI. We consider the SOI entrance state, with position $\mathbf{r}^T=[-104551.25,-597237.74,-110314.51]$, km, and velocity $\mathbf{v}^T=[3.25,17.76,3.67]$, km/s, at the time $t_0 = 8119.84$ MJD2000, in the Ecliptic J2000 reference frame centered on Venus. We simulate the flyby motion using an ephemeris model that includes the Solar System planets, the Moon and general relativity effects up to the "critical" point where $\cos\psi=-1/\sqrt{3}$. Then, we compare in Fig.~\ref{f:ThirdBodyError} the consequent forward simulations, starting from this critical point and up to the SOI exit, of the Lie transform-based approach, the numerical simulation of the second order term only of the third body potential, the Keplerian integration, and the full third body effect (not expanded).
\par

In clear contrast with the improvements obtained for the main problem case, we observe from Fig.~\ref{f:ThirdBodyError} the divergent behavior of the analytical perturbation solution. Setting the boundary condition as the Keplerian solution where the perturbing potential vanishes, in analogy with Hori's approach, the so obtained analytical solution provides even worse predictions than the Keplerian only case, without converging toward the numerically simulated trajectory for the same physical model. Other types of boundary conditions should be explored in dedicated works, for instance the magnitude difference between the perturbation and the Keplerian term in specific regions of the hyperbolic trajectory.

\section{Conclusions}

Most main problem intermediaries in the theory of artificial satellites were originally computed in polar variables, and hence are not constrained to the case of bounded orbits. Therefore, the improvements they provide in their \emph{common} realization over the Keplerian approximation encompass also the case of unbounded motion, and are clearly leveraged in their \emph{naturalized} versions. However, in the latter case, the integration ``constant'' on which the perturbation theory depends upon can no longer remain arbitrary. On the contrary, it must be determined in such a way that the perturbation solution fulfills the boundary conditions at infinity derived from the Keplerian hyperbola.
\par

In particular, we have shown that Deprit's radial intermediary provides an efficient alternative to the Keplerian approximation commonly used in flyby design. Moreover, because the appearance of secular terms is not at all of concern for hyperbolic-type motion, higher order refinements of Deprit's quasi-Keplerian approach are readily computed. Because of that, the analytical computation of unbounded perturbed Keplerian motion can take full advantage of the convergence to the true solution provided by the computation of consecutive higher orders of the perturbation approach, in which case the accuracy of the analytical solution is notably increased. 

%\section*{Compliance with ethical standards}
%\paragraph{Conflict of interest} The author declares that he has no conflict of interest.

\section*{Acknowledgements}

The research has received funding from the European Research Council under the European Union's Horizon 2020 research and innovation program (grant agreement 679086 - COMPASS). ML also acknowledges partial support from the Spanish State Research Agency and the European Regional Development Fund (Projects PID2020-112576GB-C22 and PID2021-123219OB-I00, AEI/ERDF, EU). 
%{\color{red}
%This research was motivated by discussions with C. Colombo and A. Masat, Politecnico di Milano, on the flyby dynamics and Hori's solution.
%This work is part of the dissertation that Alessandro Masat will submit to Politecnico di Milano in order to obtain the doctoral degree.
%}

\end{document}